\documentclass[reqno]{amsart}
\usepackage{graphicx,hyperref,amsthm,amssymb}
\usepackage{booktabs}
\usepackage{xcolor}
\usepackage{enumitem}
\newtheorem{theorem}{Theorem}[section]

\newtheorem{corollary}{Corollary}[section]
\theoremstyle{remark}
\newtheorem{remark}{Remark}[section]
\theoremstyle{definition}
\newtheorem{example}{Example}[section]

\title[Optimal scheduled learning rate for randomized Kaczmarz]{An optimal scheduled learning rate for a randomized Kaczmarz algorithm}

\author[N.~F.~Marshall]{Nicholas F. Marshall}
\address{Department of Mathematics, Oregon State University}
\email{marsnich@oregonstate.edu}

\author[O.~Mickelin]{Oscar Mickelin}
\address{Program in Applied and Computational Mathematics, Princeton University}
\email{hm6655@princeton.edu}

\thanks{N.F.M. was supported in part by NSF DMS-1903015.}

\begin{document}

\begin{abstract}
We study how the learning rate affects the performance of
a relaxed randomized Kaczmarz algorithm for solving $A x \approx b + \varepsilon$,
where $A x =b$ is a consistent linear system and $\varepsilon$ has independent mean zero random entries. 
We derive a learning rate schedule which optimizes a bound on the expected error that is sharp in certain cases;
in contrast to the exponential convergence of the standard randomized Kaczmarz algorithm, our optimized bound 
involves the reciprocal of the Lambert-$W$ function of an exponential.
\end{abstract}

\maketitle

\section{Introduction and main result}

\subsection{Introduction} \label{intro}
\begin{figure}[ht!]
\centering
\begin{tabular}{cc}
\includegraphics[width=.4\textwidth]{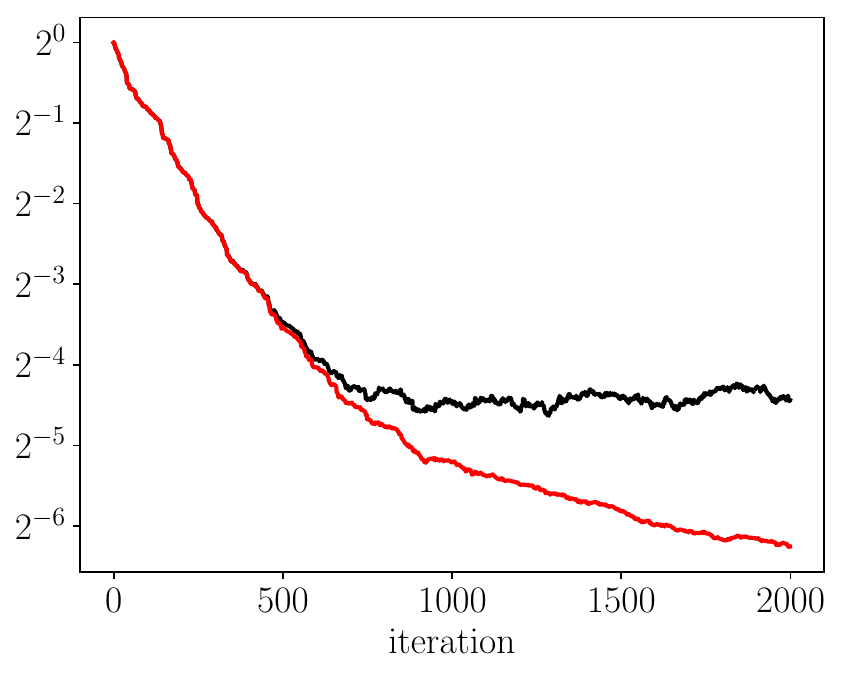} &
\includegraphics[width=.4\textwidth]{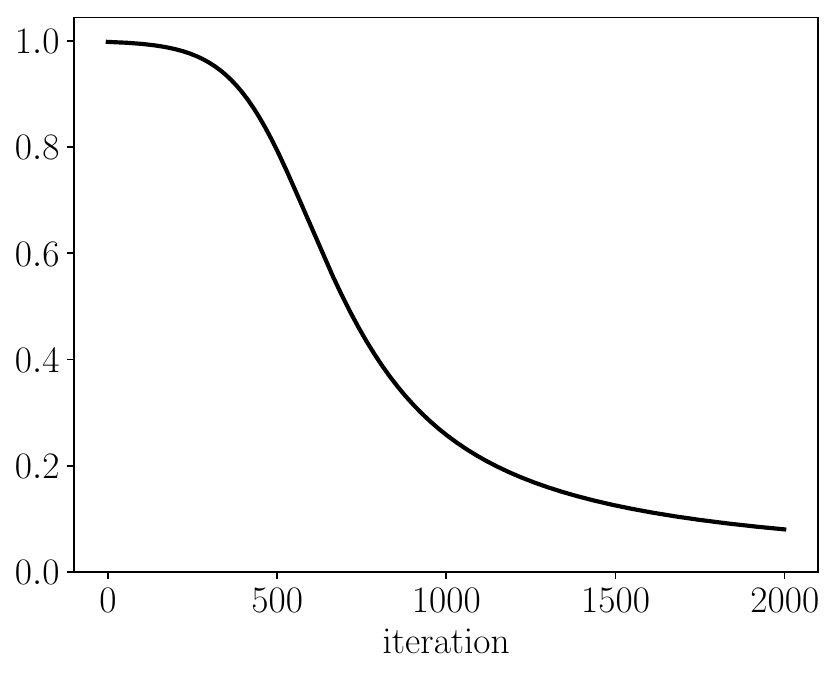}
\end{tabular}
\caption{Left: relative error when using optimal learning rate
$\alpha_k$ (red) and naive learning rate $\alpha_k=1$ (black).
Right: plot of optimal learning rate $\alpha_k$. See Example \ref{example} for details.}
\label{fig01}
\end{figure}

Let $A$ be an $m \times n$ matrix and
$A x = b$ be a consistent linear system of equations. Suppose that $\tilde{b}$
is a corrupted version of $b$ defined by
\begin{equation} \label{noise}
\tilde{b} = b + \varepsilon,
\end{equation}
where $\varepsilon$ has independent mean zero random entries. Given an initial vector $x_0$, we consider the relaxed Kaczmarz 
algorithm 
\begin{equation} \label{eq1}
x_{k+1} = x_k + \alpha_k \frac{\tilde{b}_{i_k}  - \langle a_{i_k}, x_k \rangle}{\|a_{i_k}\|^2} a_{i_k},
\end{equation}
where $\alpha_k$ is the learning rate (or relaxation parameter), $a_{i}$
is the $i$-th row of $A$, { $\tilde{b}_{i}$ is the $i$-th element of $\tilde{b}$, $i_k$ is the row index for iteration $k$, $\langle \cdot, \cdot \rangle$ is the $\ell^2$-inner product, }and
$\|\cdot\|$ is the $\ell^2$-norm. When the rows $a_{i_k}$ are chosen
randomly, \eqref{eq1} is an instance of stochastic gradient
descent, see \cite{Needell2015}, whose performance in practice depends
on the definition of the learning rate, see \cite{Goodfellow2016}. { Moreover, \eqref{eq1} can also be viewed as an instance of a stochastic Newton method, see \cite{Chung2020}.}
 In this paper, we derive a scheduled learning rate for a randomized
Kaczmarz algorithm, which optimizes a bound on the expected error; our main result
proves an associated convergence result,
see Theorem
\ref{thm1} and Figure \ref{fig01}. 

\subsection{Background}
The Kaczmarz algorithm dates back to the 1937 paper by Kaczmarz
\cite{Kaczmarz1937} who considered the
iteration \eqref{eq1} for the case $\alpha_k=1$. The algorithm was subsequently
studied by many authors; in particular, in 1967,
Whitney and Meany \cite{WhitneyMeany1967} established a convergence result for
the relaxed Kaczmarz algorithm: if $A x = b$ is a
consistent linear system, $i_k = k \mod m$, and $\alpha_k = \mu$ for fixed $0 <
\mu < 2$, then \eqref{eq1} converges to $x$, see 
\cite{WhitneyMeany1967}. 
In 1970 the Kaczmarz algorithm was rediscovered under
the name Algebraic Reconstruction Technique (ART) by  Gordon, Bender, and Herman
\cite{GordonBenderHerman1970}  who were interested in applications to
computational tomography (including applications to three-dimensional electron
microscopy); such applications typically use the relaxed Kaczmarz algorithm with
learning rate $0 < \alpha_k \le 1$, see \cite{Herman2009}. Methods and
heuristics for setting the learning rate $\alpha_k$ have been considered by
several authors, see the 1981 book by Censor \cite{Censor1981}; also see
\cite{Censor1983, Hanke1990, Hanke1990b}. 

More recently, in 2009, Strohmer and Vershynin \cite{StrohmerVershynin2009}
established the first proof of a convergence rate for a Kaczmarz algorithm that
applies to general matrices; in particular, given a consistent linear system
$Ax = b$, they consider the iteration \eqref{eq1}  with $\alpha_k =1$. Under
the assumption that the row index $i_k$ at iteration $k$ is chosen randomly with
probability proportional to $\|a_{i_k}\|^2$  they prove that
\begin{equation} \label{eqsv}
\mathbb{E} \|x_k - x\|^2 \le (1 - \eta)^{k} \|x - x_0\|^2,
\end{equation}
where $\eta := \kappa(A)^{-2}$ {and $\mathbb{E}$ is the expected value operator}; here, $\kappa(A)$ is a condition number for the
matrix $A$  defined by $\kappa(A) := \|A\|_F\|A^{-1}\|$, 
{
where $A^{-1}$ is the left inverse of $A$,
 $\|A^{-1}\|$ is the operator norm of $A^{-1}$, and
 $\|A\|_F$ is the Frobenius norm of $A$.}
We remark that the convergence rate in \eqref{eqsv} is referred to as exponential
convergence in \cite{StrohmerVershynin2009}, while in the field of numerical
analysis (where it is typical to think about error on a logarithmic scale) it
is referred to as linear convergence.

The result of \cite{StrohmerVershynin2009} was subsequently extended by Needell
\cite{Needell2010} who considered the case of a noisy linear system: instead of having
access to the right hand side of the consistent linear system $A x = b$, we 
are given $\tilde{b} = b + \gamma$, where the entries of $\gamma$ 
satisfy $|\gamma_i| \le \delta \|a_i\|$ but are otherwise arbitrary.
Under these assumptions \cite{Needell2010} proves that the
iteration \eqref{eq1} with $\alpha_k = 1$ satisfies
\begin{equation} \label{needell}
\mathbb{E} \|x_k - x\|^2 \le (1 - \eta)^{k} \|x - x_0\|^2 + \frac{\delta^2}{\eta},
\end{equation}
that is, we converge {in expectation} until we reach some ball of radius $\delta/\sqrt{\eta}$ around the solution and then no more. {Recall that} $A^{-1}$ is the left inverse of $A$, and observe that
\begin{equation} \label{esttt}
\|x - A^{-1}(b + \gamma)\| = \|A^{-1} \gamma \| \le \|A^{-1} \| \|\gamma\| \le
\delta \|A^{-1} \| \|A\|_F = \frac{\delta}{\sqrt{\eta}}.
\end{equation}
Moreover, if $\gamma$ is a scalar multiple of the left singular vector of $A$
associated with the smallest singular value, and $|\gamma_i| = \delta \|a_i\|$,
then \eqref{esttt} holds with equality (such examples are easy to manufacture).
Thus, \eqref{needell} is optimal when $\gamma$ is arbitrary. 
In this paper, we consider the case where $\tilde{b} = b + \varepsilon$,
where $\varepsilon$ has independent mean zero random entries: our main
result shows that in this case, we break through the convergence horizon
of \eqref{needell} by using an optimized learning rate and many equations with
independent noise, see Theorem \ref{thm1} for a precise statement.

{
\begin{remark}[Breaking through convergence horizon]
To be clear, when we say that our method breaks through the convergence horizon of \eqref{needell}, 
we mean that we are able to achieve expected error $\mathbb{E} \|x_k - x\|^2$ less than $\delta^2/\eta$. We achieve this by 
considering the model \eqref{eq1} where equations have independent noise and by using an optimal learning rate. 
Our main result establishes a bound on the expected error that decreases to zero as the number of available equations increases to infinity, which is the case in applications  involving streaming, such as computational tomography.
Note that under the model \eqref{eq1}, determining the solution exactly requires an infinite number of equations; indeed, this is clearly the case even if the system of equations only has a single unknown.
When a finite number of equations are available, we derive a precise bound on the expected error as a function of the number of iterations; Figure \ref{fig01} demonstrates how the optimal learning rate improves the error compared to the standard Kaczmarz algorithm in a finite number of iterations.
\end{remark}
}

\subsection{Related work}
Modifications and extensions of the randomized Kaczmarz algorithm have been considered by many authors, see
\cite{
Cai2012,
eldar2011acceleration,
Haddock2020,
Liu2015,
Ma2015,
Necoara2019,
Needell2015b,
needell2014paved,
Petra2016,
slagel2019sampled,
Steinerberger2022,
Tan2018,
Zouzias2013}. 
We note that
Cai, Zhao, and Tang \cite{Cai2012} previously considered a relaxed Kaczmarz
algorithm, but their analysis focuses on the case of a consistent linear system. 
More recently, Haddock, Needell, Rebrova, and Swartworth \cite{Haddock2020}
considered the case of a consistent linear system corrupted by sparse noise; their
main result proves convergence for a class of matrices by using an adaptive
learning rate, which roughly speaking, attempts to avoid projecting onto
corrupted equations. The result was subsequently generalized by
Steinerberger \cite{Steinerberger2022}. Our results are complementary to
the results of \cite{Haddock2020,Steinerberger2022}:  we allow for corruption of all elements of $b$, but assume that 
corruptions are independent symmetric random variables; we show 
that, by using an optimal learning rate schedule, we can recover the solution to any accuracy if we have access to a sufficient 
number of equations with independent noise.

The randomized Kaczmarz algorithm  of \cite{StrohmerVershynin2009}
can also be viewed as an instance of other machine learning methods. 
In particular, it can be viewed as an instance of coordinate descent, see
\cite{Wright2015}, or as an instance of stochastic gradient descent, see
\cite{Needell2015},
{
or an instance of the stochastic Newton method, see \cite{Chung2020, slagel2019row}.}
Part of our motivation for studying learning rate schedules
for the randomized Kaczmarz algorithm is that the randomized Kaczmarz
algorithm provides a simple model where we can start to develop a complete theoretical
understanding of the precise benefits of learning rates.
The learning rate is of essential importance in machine learning; in particular,
in deep learning: ``The learning rate is perhaps the most important hyperparameter.... the effective capacity of the model is highest when the learning
rate is correct for the optimization problem'',  \cite[pp. 429]{Goodfellow2016}.

 In this paper, we derive a scheduled learning rate (depending on two hyperparameters) for a randomized
Kaczmarz algorithm, which optimizes a bound on the expected error, and prove an associated
convergence result. Here, the word scheduled refers to the fact that the
learning rate is determined a priori as a function of the
iteration number and possibly hyperparameters (that is, the
rate is non-adaptive). See \cite{Robbins1951,Xu2011} for some general results
about learning rate schedules. 

\subsection{Summary of main contributions}
Given a consistent linear system of equations
$A x = b$, we study the problem of recovering the solution $x$ from $A$ and a corrupted
right hand side $\tilde{b}$, where $\tilde{b} = b +
\varepsilon$, and $\varepsilon$ has independent mean zero entries with bounded variance.
We show the following:
\begin{itemize}
\item In contrast to the case of a consistent linear system 
 (or a system with adversarial noise) changing the learning rate $\alpha_k$ is
advantageous. 
\item We derive a scheduled learning rate that optimizes a bound on the expected error.
\item In contrast to previous works related to the randomized Kaczmarz
algorithm that exhibit exponential convergence,
our optimized error bound involves 
the reciprocal of the Lambert-$W$ function of an exponential.
\item In the limit as the number of iterations $k$ tends to infinity, the
optimal learning rate $\alpha_k$ converges to time-based
decay  $1/(1 + \eta k)$, which
is a classic learning rate schedule. 
\item Our analysis has the potential to be combined with
modifications of the randomized Kaczmarz algorithm 
to address other related problems, {such as block-Kaczmarz methods or using the model \eqref{noise} and \eqref{eq1} to study adaptive learning rates.}
\end{itemize}

\subsection{Main result} \label{mainresult}
Let $A$ be an $m \times n$ matrix, $A x = b$ be a consistent linear system of
equations, and $a_i$  denote the $i$-th row of $A$. Suppose that $\tilde{b}$ is defined by
$$
\tilde{b} = b + \varepsilon,
$$
where $\varepsilon = (\varepsilon_1,\ldots,\varepsilon_m)$ are independent
random variables such that $\varepsilon_i$ has mean $0$ and variance
$\sigma^2 \|a_i\|^2$. 
{ This assumption about the variance can be interpreted as assuming that the data has a common signal-to-noise ratio.
}
Given $x_0$, define
\begin{equation} \label{itr}
x_{k+1} = x_k + \alpha_k \frac{\tilde{b}_{i_k}  - \langle a_{i_k}, x_k \rangle}{\|a_{i_k}\|^2} a_{i_k},
\end{equation}
where $\alpha_k$ denotes the learning rate parameter; assume that $i_k$
is chosen from $\{1,\ldots,m\} \setminus
\{i_0,\ldots,i_{k-1}\}$ with probability proportional to
$\|a_{i_k}\|^2$. Let $A_k$ denote the $(m-k) \times n$ matrix formed by
deleting rows $\{i_0,\ldots,i_{k-1}\}$ from $A$. 

\begin{theorem} \label{thm1}
In addition to the assumptions stated in \S \ref{mainresult} { above}, assume that
\begin{equation} \label{etacond}
\kappa(A_k)^{-2} \ge \eta, \quad \text{for} \quad k \le N{-1},
\end{equation}
{ for some $N \le m$.}
Given $\eta$ and $\|x - x_0\|^2/\sigma^2$ set,
$\beta_0 := \|x - x_0\|^2/\sigma^2$  and define $\alpha_k$ recursively by
\begin{equation} \label{optlearn}
\alpha_k = \frac{ \eta \beta_k  }{ \eta \beta_k  + 1 },
\quad \text{and} \quad
\beta_{k+1} = \beta_k\left(1 - \eta \alpha_k\right), \quad \text{for} \quad k=0,1,\ldots.
\end{equation}
If $x_k$ is defined iteratively by \eqref{itr}, then
\begin{equation} \label{result}
\mathbb{E} \|x_k - x\|^2 \le \sigma^2 \beta_k \le f(k), \quad \text{for} \quad k \le N,
\end{equation}
where
\begin{equation} \label{deffc}
f(k) :=
%%%% \frac{\sigma^2}{\eta W(e^{\eta k+c})},
\frac{\sigma^2}{\eta W(c e^{\eta k}{)}},
\quad \text{for} \quad
c := \frac{\sigma^2 }{\eta \|x - x_0\|^2} \exp \left( \frac{\sigma^2 }{\eta \|x - x_0\|^2} \right),
\end{equation}
where $W$ denotes the Lambert-$W$ function
(the inverse of the function $x \mapsto
xe^x$).
\end{theorem}

The proof of Theorem \ref{thm1} is given in \S \ref{proofmain}. 
In the following, we illustrate Theorem \ref{thm1} 
with a mix of remarks, examples, and corollaries.

\begin{remark}[Interpreting Theorem \ref{thm1}] \label{rmk1} We summarize a few key observations that aid in interpreting Theorem \ref{thm1}:
\begin{enumerate}[label={(\roman*)}]

\item \label{rmk1b1} Asymptotic behavior of $f(k)$: in the limit as the noise goes to zero
\begin{equation} \label{nonoise}
f(k) \sim e^{-\eta k} \|x - x_0\|^2, \quad \text{as} \quad \sigma \rightarrow 0,
\end{equation}
and in the limit as the number of iterations go to infinity
\begin{equation} \label{allnoise}
f(k) \sim \frac{\sigma^2}{\eta^2 k}, \quad \text{as} \quad k \rightarrow \infty,
\end{equation}
see Corollary \ref{cor1} and \ref{cor2}, respectively, for more precise
statements. In particular, this implies that $\mathbb{E}\|x-x_k\|^2$ tends to
$0$ as $k \rightarrow \infty$ { (letting $k$ tend to infinity
requires letting $N$ tend to infinity, which requires an increasing number of equations $m$ with independent noise}, see bullet
point \eqref{three} below).
\item  {
Informally speaking, the assumption \eqref{etacond} says that the submatrix of the remaining rows remains well-conditioned as we run the algorithm. 
}
 If $A$ is an $m \times n$ random matrix with i.i.d. rows, then
\eqref{etacond} can be replaced by a condition on the distribution of
$A$, which is independent of the iteration
number $k$. In this case, the result
\eqref{result} holds for all $k \le m$, see Corollary \ref{cor3}.
\item \label{three} We emphasize that $i_k$ is chosen without replacement, see \S
\ref{mainresult}, so the maximum number of iterations is at most the
number of rows $m$. In practice, the algorithm can be run
with restarts: after a complete pass over the data we use the final iterate to
redefine $x_0$, and restart the algorithm (potentially using different hyperparameters to define the learning rate). In the context of machine learning, 
the statement of Theorem \ref{thm1} applies to one epoch,  
see \S \ref{discuss}.
%%%\item  If $A$ is an $m \times n$ random matrix
%%%with i.i.d. rows, and $a$ denotes an arbitrary row of $A$, then \eqref{etacond} can be replaced by the condition
%%%$$
%%%\mathbb{E}  \left|\left\langle z, \frac{a}{\|a\|} \right\rangle\right|^2 \ge \eta \|z\|^2 ,
%%%\quad \forall z \in \mathbb{R}^n ,
%%%$$
%%%and the result \eqref{result} of Theorem \ref{thm1}  holds for all $k \le m$ if we loop over the rows of $A$ in order $i_k = k$,
%%%see Corollary \ref{cor3}.

\item The error bound $f(k)$ is sharp
in some cases; in these cases the learning rate is optimal, see Corollary \ref{cor4}.
We demonstrate this corollary numerically in Figure \ref{fig02}.
\end{enumerate}
\end{remark}

\begin{example}[Numerical illustration of Theorem \ref{thm1}] \label{example}
Let $A$ be an $m \times n$ matrix with $s$ nonzero entries in
each row. Assume these nonzero entries are independent vectors drawn
uniformly at random from the unit sphere: $S^{s-1} = \{ v
\in \mathbb{R}^s : \|v\|=1 \}$. Let $x$ be an $n$-dimensional vector with
independent standard normal entries, and set $b := A x$. Let
$\tilde{b} = b + \varepsilon$, where $\varepsilon$ is an $m$-dimensional vector
with independent mean $0$, variance $\sigma^2$ normal entries. We run the
relaxed Kaczmarz  algorithm \eqref{eq1} using the learning rate
\eqref{optlearn} of Theorem \ref{thm1}. In particular, we set
$$
m =2000, \quad n =100, \quad s=10, \quad \sigma = .05, \quad \text{and} \quad x_0 = \vec{0}.
$$
Using the estimates  $\|x - x_0\|^2 = n$ and $\eta = 1/100$ we define $\alpha_k$ by \eqref{optlearn}. This choice of $\eta$ is justified by  Corollary \ref{cor4} below in combination with the fact that the rows of $A$ are isotropic when scaled by $\sqrt{n}$, see \cite[\S 3.2.3, \S 3.3.1]{vershynin2018high}.
We plot the numerical relative error $\|x-x_k\|/\|x\|$ together with the bound on the expected relative error $\sqrt{f(k)}/\|x\|$ in Figure \ref{fig02}.
Furthermore, to provide intuition about how $f(k)$ varies with $\sigma$, we plot $\sqrt{f(k)}/\|x\|$ for various values of
$\sigma$ in Figure \ref{fig02}, keeping
other parameters fixed.

{At first, when
$\|x - x_k\|^2 \gg \sigma^2$, (roughly iterations 1 to 500, see Figure \ref{fig01})  the error decreases linearly in the
logarithmic scale of the figure illustrating  the asymptotic rate 
\eqref{nonoise}.
For large $k$, when $\|x - x_k\|^2 \ll \sigma^2$ (after roughly 1000 iterations, see Figure \ref{fig01}) the error decreases like $1/\ln k$ in the logarithmic scale of the figure, illustrating the asymptotic rate \eqref{allnoise}.
}
\begin{figure}[h!]
\centering
\begin{tabular}{cc}
\includegraphics[width=.4\textwidth]{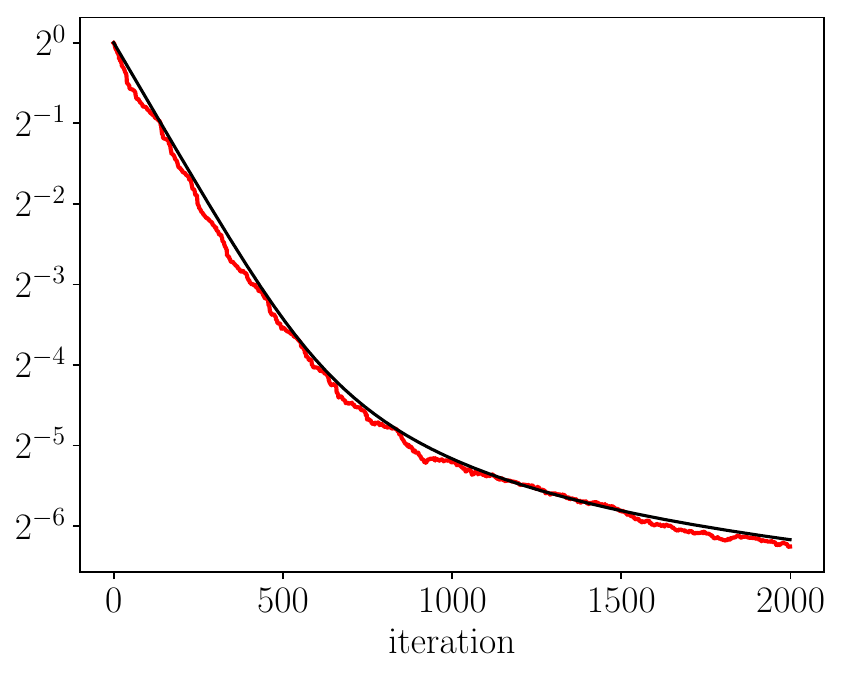}  &
\includegraphics[width=.4\textwidth]{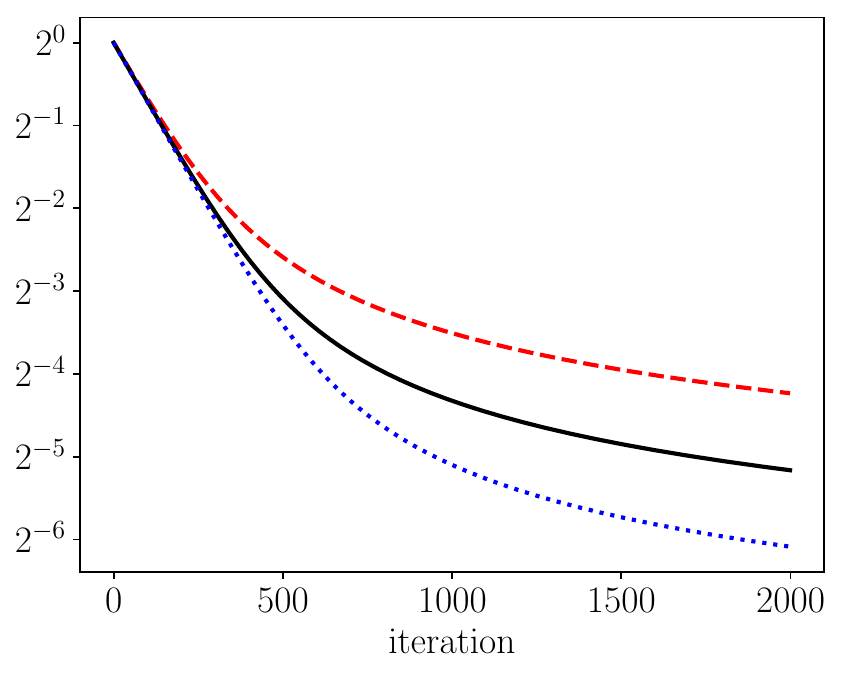} 
\end{tabular}
\caption{Left: relative error $\|x-x_k\|/\|x\|$ (red), and  function $\sqrt{f(k)}/\|x\|$ 
(black). Right: how the function $\sqrt{f(k)}/\|x\|$ changes when
$\sigma = .05$ (dotted blue), $\sigma = .1$ (solid black), and $\sigma = .2$ (dashed red).
See Example \ref{example} for details.}
 \label{fig02}
\end{figure}
\end{example}

\begin{figure}[ht!]
\centering
\begin{tabular}{cc}
\includegraphics[width=.4\textwidth]{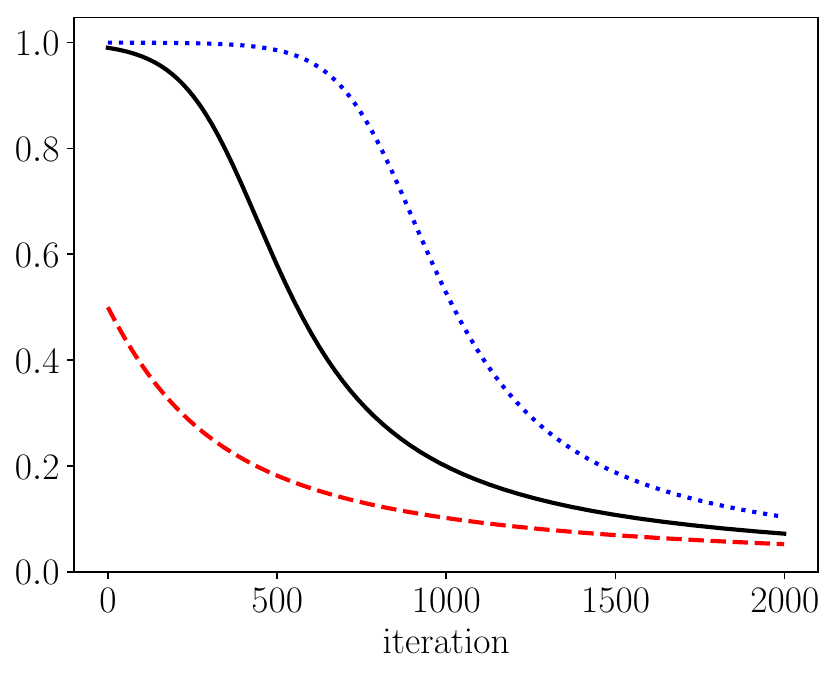} &
\includegraphics[width=.4\textwidth]{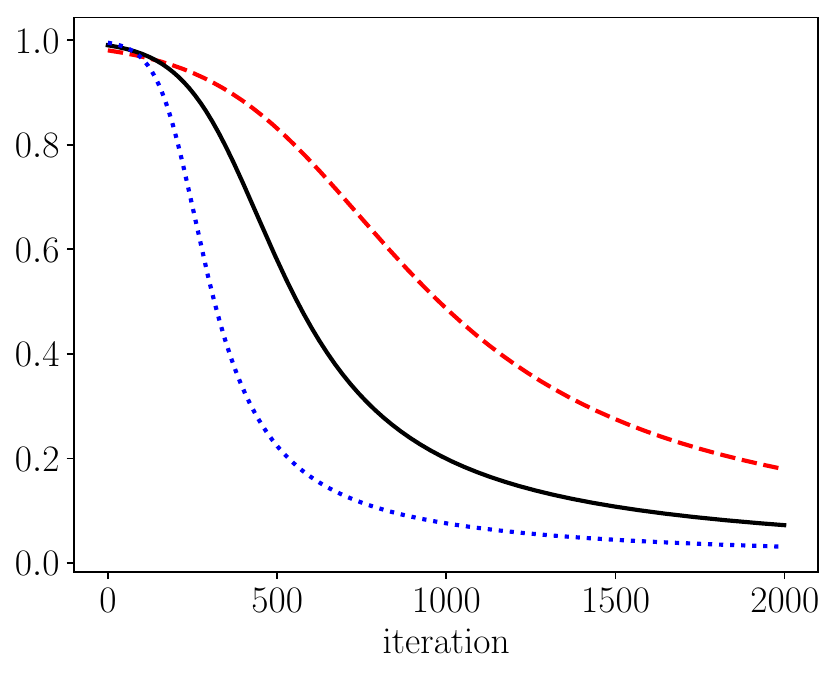}
\end{tabular}
\caption{Left: $\alpha(t)$ for $\sigma = .01$
(dotted blue), $\sigma = .1$ (solid black), and $\sigma = 1$ (dashed red).
Right: $\alpha(t)$ for $\eta = .005$ (dashed red), $\eta = .01$  (solid black),and $\eta = .02$ (dotted blue).}
\label{fig03}
\end{figure}

\begin{example}[Continuous version of learning rate $\alpha_k$] \label{numerics}
The learning rate $\alpha_k$ defined in \eqref{optlearn} optimizes the error bound of Theorem \ref{thm1}, see \S \ref{step4}.
The scheduled learning rate $\alpha_k$ depends on two parameters (assuming $x_0 = \vec{0}$):
\begin{itemize}
\item the signal-to-noise ratio $\|x\|^2/\sigma^2$, and
\item the condition number parameter $\eta$.
\end{itemize}
The result of Theorem \ref{thm1} states that the function $f(k)$ defined in \eqref{deffc}
is an upper bound for $\beta_k \sigma^2$. From the proof of Theorem \ref{thm1} it will be clear that
this upper bound is a good approximation when $\eta$ is small, see \S \ref{step5}. In this case, it is
illuminating to consider a continuous version of the optimal scheduled learning rate of Theorem \ref{thm1}. 
In particular, we define
\begin{equation} \label{alphacont}
\alpha(t) = \frac{\eta f(t)}{\eta f(t) + \sigma^2},
\end{equation}
where $f(t)$ is defined by \eqref{deffc}. We plot the function $\alpha(t)$ for three different levels of noise: $\sigma = .01$, $\sigma = .1$, and $\sigma = 1$, while keeping the other parameters fixed (and set by the values in Example \ref{example}),
 see Figure \ref{fig03}. 
\end{example}

 We also plot the function $\alpha(t)$ for three different values of the condition number parameter: 
 for $\eta = .005$, $\eta = .01$ and $\eta = .02$,
 while keeping the other parameters fixed (and set by the values in Example \ref{example}),
 see Figure \ref{fig03}. 
 
\begin{remark}[Asymptotics of $\alpha(t)$]
The continuous version of the learning rate \eqref{alphacont} has two distinct asymptotic regimes similar to
 Remark \ref{rmk1} \eqref{rmk1b1}. In particular, 
in the limit as the noise goes to zero we have
\begin{equation} \label{alphasigmasmall}
\alpha(t) \sim 1, \quad \text{as} \quad \sigma \rightarrow 0,
\end{equation}
and in the limit as the number of iterations goes to infinity we have
\begin{equation} \label{alphakbig}
{\alpha}(t) \sim \frac{1}{1 + \eta {t}}, \quad \text{as} \quad {t} \rightarrow \infty.
\end{equation}
The regimes \eqref{alphasigmasmall} and \eqref{alphakbig} are illustrated by
the blue dotted line in the left and right plots of Figure
\ref{fig03}, respectively.
We note that \eqref{alphakbig} corresponds to time-based decay,
which is a popular learning rate in practice, see the learning rate schedules of \cite{Tensorflow}.
\end{remark}

{
\begin{example}[Optimal learning rate versus time-based decay]

 We compare the optimal learning rate $\alpha_k$ to the time-based decay learning rate $1/(1+\eta k)$, which is the large iteration limit of the optimal learning rate,  see \eqref{alphakbig}. We run the Kaczmarz algorithm \eqref{eq1} with these learning rates for the system described in Example \ref{example}, see Figure \ref{fig04new}.
 
\begin{figure}[h!]
\begin{tabular}{cc}
\includegraphics[width=.4\textwidth]{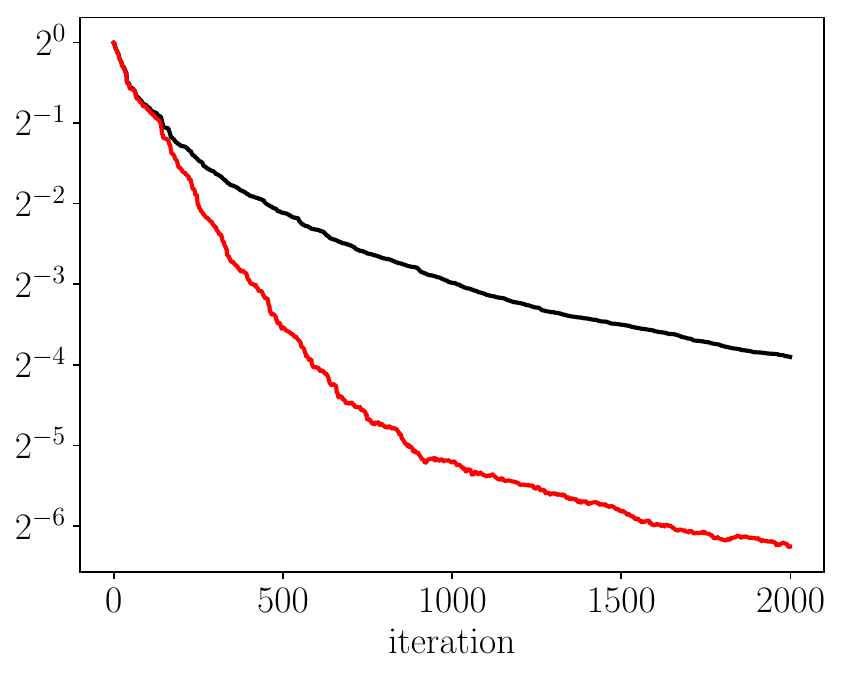} &
\includegraphics[width=.4\textwidth]{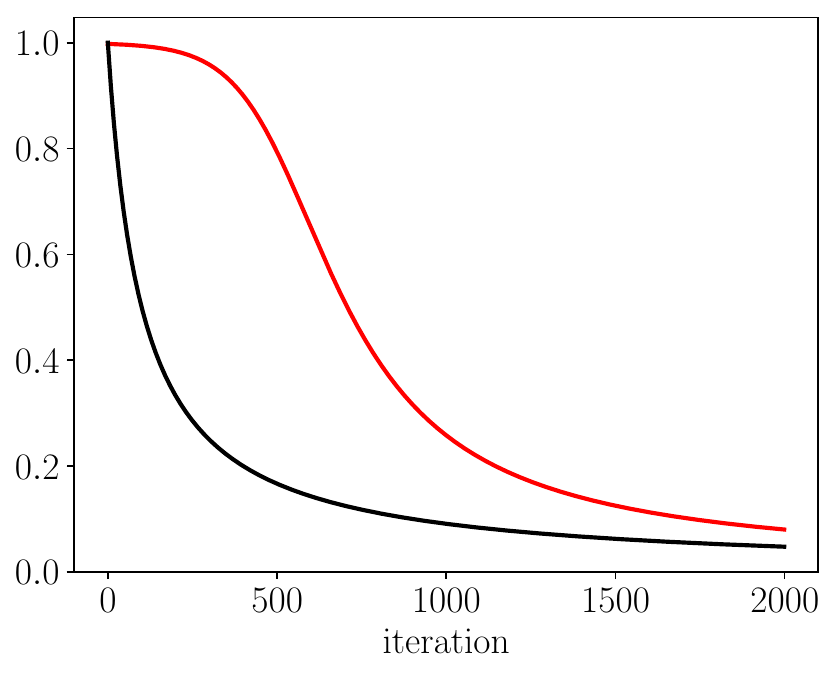}
\end{tabular}
\caption{Left: relative error when using optimal learning rate
$\alpha_k$ (red) and time-based decay $1/(1+\eta k)$ (black). Right: plot of  optimal learning rate $\alpha_k$ (red), and time-based decay (black)} \label{fig04new}
\end{figure}

The fact that time-based decay $1/(1+\eta k)$ is the large iteration limit of the optimal learning rate $\alpha_k$ is reflected in the fact that eventually (after $1500$ iterations) both errors appear to decrease at similar rates. However, initially the optimal learning rate decreases the error much faster and a gap between the two relative errors appears (and this gap will remain).
Informally speaking, the reason that using an `S'-shaped learning rate is optimal is that initially, when the error is larger than the noise, it is advantageous
to keep the learning rate close to $1$ (to make rapid initial progress), and to only decrease the learning rate once the error is smaller than the noise. 
It is instructive to compare these results to Figure \ref{fig01}, where we plot 
the errors of the optimal learning rate to the constant learning rate $1$, which is the small noise limit of the optimal learning rate, see \eqref{alphasigmasmall}. Observe that in Figure \ref{fig01}, the errors initially decrease at the same rate, but eventually the error associated with the constant learning rate stagnates. Informally speaking, the optimal learning rate can be viewed as an optimal transition between the constant learning rate $1$ and time-based decay $1/(1+\eta k)$.
\end{example}
}

{
\begin{example}[Estimating hyperparameters]\label{ex:est_hyperaparams}
As noted in Example \ref{numerics}, the scheduled learning rate $\alpha_k$ defined in Theorem \ref{thm1} depends on two hyperparameters
the signal-to-noise ratio $\|x\|^2/\sigma^2$, and the condition number parameter $\eta$.
In order to use this learning rate in practice, it is necessary to estimate these parameters. Indeed, commonly used Learning rate schedules such as the constant learning rate, time-based decay, step-based decay, and exponential decay depend on one or more parameters
\cite{Tensorflow}, and the problem of tuning these hyperparameters is important in practice \cite{Goodfellow2016}. The fact that the hyperparameters of the of the learning rate schedule $\alpha_k$ of Theorem \ref{thm1} have natural interpretations 
as the signal-to-noise-ratio and condition number parameter, respectively, provides addition intuition towards tuning these parameters.

We next show a heuristic to estimate the parameters $\eta$ and $\|x\|^2/\sigma^2$ from one additional run of the randomized Kaczmarz method. Let $x_0,\ldots,x_N$ denote the iterates resulting from running the 
randomized Kaczmarz algorithm \eqref{eq1} with constant learning rate $1$ for $N$ iterations. If we assume the initial error is above the noise level, and assume that eventually the error stagnates because of the noise, then we expect that, initially, $\|x_{j+1} - x\|^2 \approx (1 - \eta) \|x_j - x\|^2$, see \eqref{Eiter},
and eventually $\|x_j - x\|^2 \approx \sigma^2/\eta$, see \eqref{needell}. Using these estimates and the approximation $\|x_N\|^2 \approx \|x\|^2$ yields the heuristic estimates
for $\eta$ and $\beta_0 = \|x\|^2/\sigma^2$ 
\begin{equation}\label{eq:heuristic_parameters}
\tilde{\eta} := 1 - \frac{1}{N_0} \sum_{j=1}^{N_0} \frac{\|x_j - x_N\|^2}{\|x_{j-1} - x_N\|^2},
\quad \text{and} \quad
\tilde{\beta}_0^{-1} := \tilde{\eta} \frac{1}{N_1} \sum_{j=N-N_1}^{N-1} \frac{\|x_N - x_j\|^2}{\|x_N\|^2},
\end{equation}
for some integers $N_0$ and $N_1$. Suppose that the system described in Example \ref{example} is given, but the parameter condition number parameter $\eta$ and the signal-to-noise ratio $\|x\|^2/\sigma^2$ are unknown.
In Figure~\ref{fig05}, we show the result of running the Kaczmarz method with the estimated learning rate $\tilde{\alpha}_k$ obtained from the parameters estimated by \eqref{eq:heuristic_parameters}, using $N_0 = 100$, $N_1 = 100$ and $N=2000$.
\begin{figure}[h!]
\centering
\begin{tabular}{cc}
\includegraphics[width=.4\textwidth]{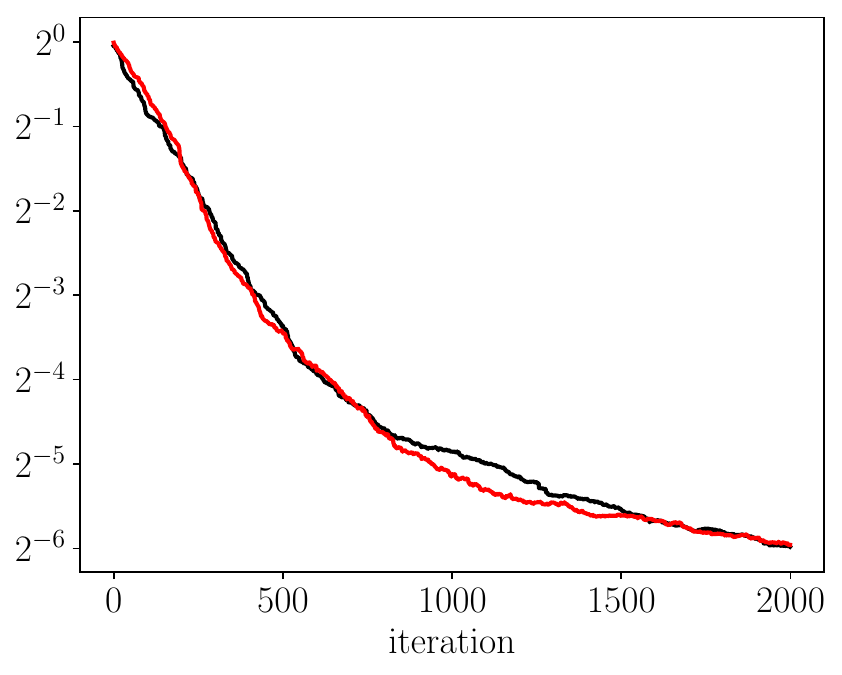} &
\includegraphics[width=.4\textwidth]{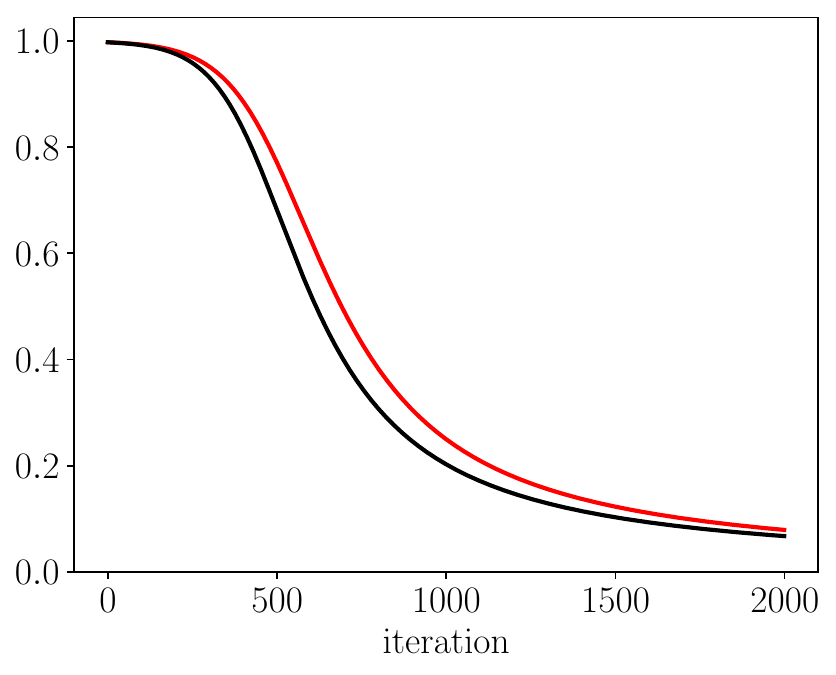}
\end{tabular}
\caption{Left: relative error when using optimal 
$\alpha_k$ (red)  and estimated $\tilde{\alpha}_k$ (black) learning rate 
Right: plot of optimal $\alpha_k$ (red) 
and estimated  $\tilde{\alpha}_k$ (black)
learning rate. }
\label{fig05}
\end{figure}

\end{example}
}

{
\begin{remark}[Normalizing rows]\label{rem:norm_rows}
Practically speaking, before the algorithm starts, the rows of the linear system can simply be normalized so that they have the same norm. We note that this can be considered as a form of preconditioning, which will change the condition number of the matrix. This row normalization simplifies the application of the randomized Kaczmarz method and can be performed as the rows are sampled in a streaming setting, as assumed in for example \cite{Haddock2020}.
\end{remark}
}

{
\begin{example}[Computational tomography example]
We lastly include a real-world example of a tomography simulation of a system of equations $A x = b$. The matrix $A$ has dimensions $4500 \times 225$ and corresponds to absorption along a random line through a $15\times 15$ grid. 
The matrix $A$ and the vector $x$ were generated from the Matlab Regularization Toolbox by P.C. Hansen \cite{Hansen2007} and $x$ was normalized to have unit norm. The right hand side $b$ of the consistent system generated by the simulation is corrupted by adding noise according to \eqref{noise}.

The numerical example combines the considerations in the preceding remarks and examples: we normalize the rows of $A$ (see Remark~\ref{rem:norm_rows}), the hyperparameters $\eta$ and $\|x\|^2/\sigma^2$ are estimated (see Remark~\ref{ex:est_hyperaparams}) and multiple epochs are used. Note that while Theorem~\ref{thm1} only applies to one epoch (see Remark~\ref{rmk1}(iii)) since we require equations with independent noise, practically speaking the result may continue to hold over multiple epochs as long as
\eqref{expectnoise} in the Proof of Theorem \ref{thm1} approximately holds.

The results are shown in Figure~\ref{fig06} for two different values of the noise. Note the similarities between Figure~\ref{fig06} and Figure~\ref{fig01}, indicating that the framework of the article applies well beyond the conditions of Theorem~\ref{thm1}.

\begin{figure}[h!]
\centering
\begin{tabular}{cc}
\includegraphics[width=.4\textwidth]{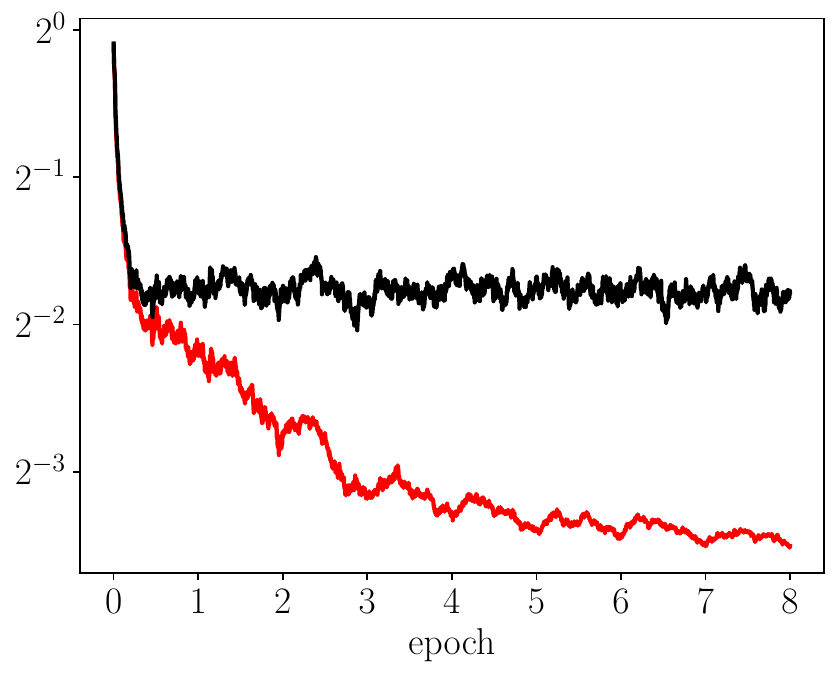} &
\includegraphics[width=.4\textwidth]{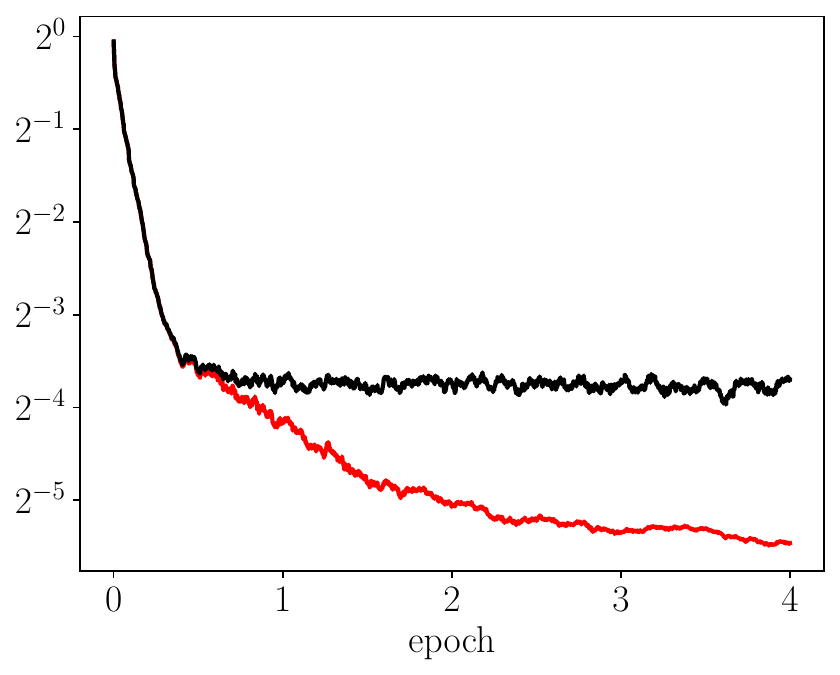}
\end{tabular}
\caption{
Relative error of constant learning rate  when using constant learning rate $1$ (black)  and estimated  learning rate (red) 
for $\sigma = 0.02$ (left).
and $\sigma = 0.005$ (right). }
\label{fig06}
\end{figure}

\end{example}
}

\subsection{Corollaries of the main result}
The following corollaries provide more intuition about the error bound function $f(k)$, and cases when the error bound of Theorem \ref{thm1} is 
sharp.
First, we consider the case where the variance of the noise $\sigma^2$ is small, and the other parameters are fixed;
in this case
we recover the convergence rate of the standard randomized Kaczmarz algorithm
with $\eta = \kappa(A)^{-2}$.

\begin{corollary}[Limit as $\sigma \rightarrow 0$] \label{cor1}
We have 
$$
f(k) = e^{-\eta k} \|x - x_0\|^2 \left(1 + \mathcal{O}\left( \frac{\sigma^2}{\|x-x_0\|^2} \right) \right), \quad \text{as} \quad \sigma \rightarrow 0,
$$
where the constant in the big-$\mathcal{O}$ notation depends on $\eta$ and $k$.
\end{corollary}

The proof of Corollary \ref{cor1} is given in \S \ref{proofcor1}. Next, we
consider the convergence rate as the number of iterations goes to infinity and
the other parameters are fixed.

\begin{corollary}[Limit as $k \rightarrow \infty$] \label{cor2}
We have
$$
f(k) = \frac{\sigma^2}{\eta^2 k} \left( 1 + \mathcal{O}\left( \frac{\ln k}{k} \right) \right),
\quad \text{as} \quad k \rightarrow \infty,
$$
where the constant in the big-$\mathcal{O}$ notation depends on $\eta$, $\|x-x_0\|^2$ and $\sigma^2$.

\end{corollary}

The proof of Corollary \ref{cor2} is given in \S \ref{proofcor2}.
Informally speaking, in combination with Theorem \ref{thm1} {and Jensen's inequality}, this corollary
says that we should expect
$$
\mathbb{E} \|x - x_k \| \lesssim \frac{\sigma}{\eta} \frac{1}{\sqrt{k}},
$$
which agrees with the intuition { (based on the central limit theorem) that using $k$ independent sources
of noise should reduce the expected error  by a factor of $1/\sqrt{k}$. }

\begin{corollary}[Matrices with i.i.d. rows] \label{cor3}
Suppose that $A$ is an $m \times n$ random matrix
with i.i.d. rows, {and let $a_0$ be a random variable generated according to this distribution}. If we run the algorithm \eqref{eq1} with $i_k=k$, and
in place of the condition \eqref{etacond} assume that
\begin{equation} \label{cond}
\mathbb{E}
\left| \left\langle z, \frac{{a_0}}{\|{a_0}\|} \right\rangle \right|^2 \ge \eta \|z\|^2,  \quad \forall z \in \mathbb{R}^n,
\end{equation}
for some fixed value $\eta > 0$, {where the expectation is over the random variable {$a_0$}.
Then the result \eqref{result} of Theorem \ref{thm1} holds for $k \le m$.}
\end{corollary}

The proof of Corollary \ref{cor3} is given in \S \ref{proofcor3}.
The condition of Corollary \ref{cor3} holds if, for example, the rows of $A$
are sampled uniformly at random from the unit sphere $S^{n-1} = \{x \in
\mathbb{R}^n : \|x\|=1\}$. Indeed, in this case 
$$
\mathbb{E}
\left| \left\langle z, a_{0}  \right\rangle \right|^2 \ge \frac{1}{n} \|z\|^2, \quad \forall z \in \mathbb{R}^n,
$$
see \cite[Lemma 3.2.3, \S 3.3.1]{vershynin2018high}. The following corollary
gives a condition under which the learning rate \eqref{optlearn}
is optimal. In particular, this corollary implies that the error bound and learning rate are optimal
for the example of matrices whose rows are sampled uniformly at random from the unit sphere discussed above.

\begin{corollary}[Case when error bound is sharp and learning rate is optimal] \label{cor4}
Assume that
\begin{equation} \label{condcond}
\mathbb{E}_{i_0,\ldots,i_{k-1}}
\left| \left\langle \frac{x_k - x,
}{\|x - x_k\|}, \frac{a_{i_k}}{\|a_{i_k}\|} \right\rangle \right|^2 = \eta.
\end{equation}
Then,
$$
\mathbb{E} \|x-x_k\|^2 = \sigma^2 \beta_k,
$$
and the learning rate $\alpha_k$ defined by \eqref{optlearn} is optimal in the sense that it minimizes the 
expected error $\mathbb{E} \|x-x_k\|^2$ over all possible choices of scheduled learning rates 
(learning rates that only depend on the iteration number $k$ and possibly hyperparameters).
Moreover, if \eqref{cond} holds with $\eta=1/n$, then it follows that \eqref{condcond} holds with equality
and hence the learning rate is optimal.
\end{corollary}
The proof of Corollary \ref{cor4} is given in \S \ref{proofcor4}.
Informally speaking, the learning rate $\alpha_k$ defined by \eqref{optlearn} is optimal whenever the 
 bound $\beta_k \sigma^2$ in Theorem \ref{thm1} is a good approximation of
the expected error $\mathbb{E} \|x-x_k\|^2$, and there is reason to expect this
is the case in practice, see the related results for consistent systems
\cite[\S 3, Theorem 3]{StrohmerVershynin2009} and
\cite[Theorem 1]{Steinerberger2021}.
For additional discussion about Theorem \ref{thm1} and its corollaries see \S
\ref{discuss}.

\section{Proof of Theorem \ref{thm1}} \label{proofmain}
{
The proof of Theorem \ref{thm1} is divided into five steps: 
\begin{itemize}
\item  Step 1 (\S\ref{step1}) We
consider the relaxed Kaczmarz algorithm for the consistent linear system $A x = b$ and prove a recursive formula for the expected error of the solution. 
This step uses the same proof strategy as 
\cite{StrohmerVershynin2009}.
\item Step 2 (\S\ref{step2}). We consider the effect of additive
random noise. More precisely, we study how the additive noise changes the analysis of Step 1. The end result is additional terms involving conditional expectations of a geometric quantity $\eta_k$ and a noise quantity $\zeta_k$, with respect to the choices of rows $i_0, \ldots , i_{k-1}$.
\item Step 3 (\S\ref{step3}). We estimate the conditional expectations of the geometric quantity $\eta_k$ and a noise quantity $\zeta_k$.
\item Step 4 (\S\ref{step4}). We optimize the learning rate $\alpha_k$ with respect to the error bound from the previous step. This results in a recurrence relation for the optimal learning rate.
\item Step 5 (\S\ref{step5}). We show that this optimized learning rate is related to a
differential equation, which can be used to establish the upper bound $f(k)$ on the expected error. 
\end{itemize}
}

\subsection{Step 1: relaxed randomized Kaczmarz for consistent systems} \label{step1}
\begin{figure}[ht!]
\centering
\includegraphics[width=.5\textwidth]{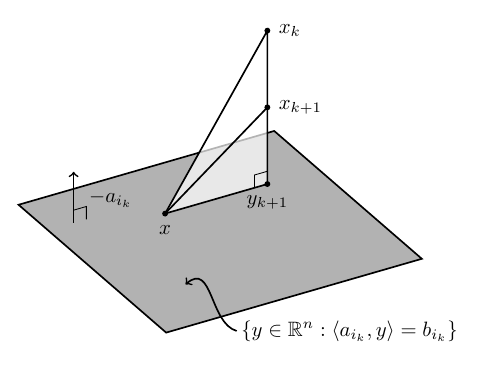}
\caption{Illustration of the points $x,y_{k+1},x_{k+1}$, and $x_{k}$.}
\label{fig04}
\end{figure}

We start by considering the consistent linear system $Ax = b$. 
Assume that $x_0$ is given and let
$$
x_{k+1} := x_k + \alpha_k \frac{b_{i_k}  - \langle a_{i_k}, x_k
\rangle}{\|a_{i_k}\|^2} a_{i_k},
$$
denote the iteration of the relaxed randomized Kaczmarz algorithm with learning rate $\alpha_k$ on the
consistent linear system $Ax = b$, and let
$$
y_{k+1} := x_k + \frac{b_{i_k}  - \langle a_{i_k}, x_k
\rangle}{\|a_{i_k}\|^2} a_{i_k},
$$
be the projection of $x_k$ on the affine hyperplane defined by the $i_k$-th
equation. The points $y_{k+1} ,x_{k+1}$, and $x_k$ lie on the
line $\{x_k + t a_{i_k} : t \in \mathbb{R} \}$, which is perpendicular 
to the affine hyperplane $\{ y \in \mathbb{R}^n : \langle a_{i_k}, y \rangle =
b_{i_k} \}$ that contains $x$ and $y_{k+1}$, 
 see the illustration in Figure \ref{fig04}.

\noindent By the Pythagorean theorem, it follows that
$$
\|x - x_{k+1}\|^2 = \|x - x_k\|^2 - \|y_{k+1} - x_k\|^2 +
\|y_{k+1} - x_{k+1}\|^2,
$$
and by definition of $y_{k+1}$ and $x_{k+1}$ we have
$$
y_{k+1} - x_{k+1}  = (1-\alpha_k) (y_{k+1} - x_{k}),
$$
see Figure \ref{fig04}. Thus
$$
\|x - x_{k+1}\|^2 = \|x - x_k\|^2 - (2\alpha_k - \alpha_k^2) \|y_{k+1} -
x_k\|^2. 
$$
Factoring out $\|x - x_k\|^2$ from the right hand side gives
$$
\|x - x_{k+1}\|^2 = \left(1 - (2\alpha_k - \alpha_k^2)  \frac{\|y_{k+1} - x_k\|^2}{\|x
- x_k\|^2} \right) \|x - x_k\|^2 .
$$
 Since $\|y_{k+1} - x_k\| = |\langle
x_k - x, a_{i_k}/\|a_{i_k}\| \rangle|$,  it follows that
$$
\|x_{k+1} - x\|^2 = \left(1 - (2\alpha_k - \alpha_k^2) \left| \left\langle \frac{x_k - x
}{\|x_k - x\|}, \frac{a_{i_k}}{\|a_{i_k}\|} \right\rangle \right|^2 \right) \|x_k - x\|^2.
$$
Taking the expectation conditional
on $\{i_0,\ldots,i_{k-1}\}$ gives
\begin{equation} \label{Eiter}
\mathbb{E}_{i_0,\ldots,i_{k-1}}
\|x_{k+1} - x\|^2 = \left(1 - (2\alpha_k - \alpha_k^2) \eta_k 
\right) \|x_k - x\|^2,
\end{equation}
where
$$
\eta_k :=
\mathbb{E}_{i_0,\ldots,i_{k-1}} \left|
\left\langle \frac{x_k - x }{\|x_k - x\|}, \frac{a_{i_k}}{\|a_{i_k}\|} \right\rangle \right|^2.
$$
We delay discussion of conditional expectation until {Step 3 (\S\ref{step3}).}

\begin{remark}[Optimal learning rate for consistent linear systems] \label{sigma0}
Observe that
$$
(1 - \eta_k) \|x_k - x\|^2 \le \left(1 - (2\alpha_k -
\alpha_k^2) \eta_k \right) \|x_k - x\|^2,
$$
with equality only when $\alpha_k=1$. It follows that, for consistent linear systems of equations, the
optimal way to define the learning rate is to set $\alpha_k = 1$ for all $k=0,1,2,\ldots$. In the following, we show that, under our noise model, defining
$\alpha_k$ as a specific decreasing function of $k$ is advantageous.
\end{remark}

\subsection{Step 2:  relaxed randomized Kaczmarz for systems with noise} \label{step2}
In this section, we redefine $x_k$ and $y_k$ for the case where 
the right hand
side of the consistent linear system $A x = b$ is corrupted by additive random
noise: $\tilde{b} = b + \varepsilon$. Let
$$
x_{k+1} : = x_k + \alpha_k \frac{\tilde{b}_{i_k}  - \langle a_{i_k},
x_k \rangle}{\|a_{i_k}\|^2} a_{i_k} =
x_k + \alpha_k \frac{b_{i_k} + \varepsilon_{i_k}  - \langle a_{i_k},
x_k \rangle}{\|a_{i_k}\|^2} a_{i_k},
$$
be the iteration of the relaxed randomized Kaczmarz algorithm using $\tilde{b}$,
and 
$$
y_{k+1} := x_k + \frac{b_{i_k}  - \langle a_{i_k}, x_k
\rangle}{\|a_{i_k}\|^2} a_{i_k},
$$
be the projection of $x_k$ on the affine hyperplane defined by the uncorrupted $i_k$-th
equation.  Note that both $x_{k+1}$ and $y_{k+1}$ differ from the previous
section:  $x_{k+1}$ is corrupted by the noise term $\varepsilon_{i_k}$ and
$y_{k+1}$ is the projection of the previously corrupted iterate $x_k$ onto the hyperplane defined by the uncorrupted equation.
However, the following expansion still holds:
\begin{equation} \label{expando}
\|x - x_{k+1}\|^2 = \|x - x_k\|^2 - \|y_{k+1} - x_k\|^2 +
\|y_{k+1} - x_{k+1}\|^2.
\end{equation}
Indeed, $x_{k+1},y_{k+1},$ and $x_k$ are still contained on the line  $\{x_k + t
a_{i_k} : t \in \mathbb{R} \}$, which is perpendicular  
to the affine hyperplane $\{ y \in \mathbb{R}^n : \langle a_{i_k}, y \rangle =
b_{i_k} \}$ that contains $x$ and $y_{k+1}$, see  Figure \ref{fig04}.
By the definition of $y_{k+1}$ and $x_{k+1}$ we have
$$
\|y_{k+1} - x_{k+1}\|^2  = \left\| 
(1-\alpha_k) (y_{k+1} - x_{k}) -  \alpha_k \frac{\varepsilon_{i_k}}{\|a_{i_k}\|} \frac{a_{i_k}}{\|a_{i_k}\|} \right\|^2.
$$
Expanding the right hand side gives 
\begin{equation} \label{expande}
\|y_{k+1} - x_{k+1}\|^2  = 
(1 - \alpha_k)^2 \| y_{k+1} - x_{k}\|^2  + Z_k,
\end{equation}
where
$$
Z_k := -2  \alpha_k(1 - \alpha_k)\frac{\varepsilon_{i_k}}{\|a_{i_k}\|}
 \left\langle y_{k+1} - x_k , \frac{a_{i_k}}{\|a_{i_k}\|} \right\rangle+
\alpha_k^2 \frac{\varepsilon_{i_k}^2}{\|a_{i_k}\|^2}.
$$
By using the fact that
$$
 \left\langle y_{k+1} - x_k , \frac{a_{i_k}}{\|a_{i_k}\|} \right\rangle = \frac{b_{i_k} - \langle a_{i_k}, x_k \rangle}{\|a_{i_k}\|},
$$
we can rewrite $Z_k$ as
$$
Z_k = -2  \alpha_k(1 - \alpha_k)\frac{\varepsilon_{i_k}}{\|a_{i_k}\|}
\frac{b_{i_k} - \langle a_{i_k}, x_k \rangle}{\|a_{i_k}\|} +
\alpha_k^2 \frac{\varepsilon_{i_k}^2}{\|a_{i_k}\|^2}.
$$
Combining \eqref{expando} and \eqref{expande} gives
$$
\|x - x_{k+1}\|^2 = \|x - x_k\|^2 - 
(2\alpha_k -\alpha_k^2) \|y_{k+1} - x_k\|^2 
+ Z_k.
$$
As in the analysis of the consistent linear system in Step 1 (\S \ref{step1}) above, we factor out
$\|x - x_k\|^2$, use the fact that $\|y_{k+1} - x_k\| = \langle
x_k - x, a_{i_k}/\|a_{i_k}\| \rangle$, and take the expectation conditional on
$\{i_0,\ldots,i_{k-1}\}$ to conclude that
\begin{equation} \label{expectn}
\mathbb{E}_{i_0,\ldots,i_{k-1}}
\|x_{k+1} - x\|^2 = \big(1 - (2 \alpha_k - \alpha_k^2) \eta_k
 \big) \|x_k - x\|^2
+ 
\zeta_k,
\end{equation}
where
$$
\eta_k := 
\mathbb{E}_{i_0,\ldots,i_{k-1}}
\left| \left\langle \frac{x_k - x
}{\|x_k - x\|}, \frac{a_{i_k}}{\|a_{i_k}\|} \right\rangle \right|^2,
\quad \text{and} \quad
\zeta_k = 
\mathbb{E}_{i_0,\ldots,i_{k-1}} Z_k.
$$
In the following section, we discuss the terms $\zeta_k$ and $\eta_k$.

\subsection{Step 3: estimating conditional expectations}\label{step3}
In this section, we discuss estimating the conditional expectations $\eta_k$ and
$\zeta_k$. First, we discuss $\zeta_k$, which has two terms (a linear term and quadratic
term with respect to $\varepsilon_{i_k}$). In particular, we have 
\begin{equation}\label{eq:def_zeta_k_new}
  \zeta_k = -2  \alpha_k(1 - \alpha_k)
\mathbb{E}_{i_0,\ldots,i_{k-1}} \frac{\varepsilon_{i_k}}{\|a_{i_k}\|}
\frac{b_{i_k} - \langle a_{i_k}, x_k \rangle}{\|a_{i_k}\|} +
\alpha_k^2 \mathbb{E}_{i_0,\ldots,i_{k-1}}  \frac{\varepsilon_{i_k}^2}{\|a_{i_k}\|^2} .
\end{equation}
Recall that $\varepsilon_1,\ldots,\varepsilon_m$ are independent random
variables such that $\varepsilon_{i_{k}}$ has mean zero and variance $\sigma^2 \|a_{i_k}\|^2$.
Since we assume that $i_k$ is chosen from $\{1,\ldots,m\}
\setminus \{i_0,\ldots,i_{k-1}\}$ (that is, they are drawn without replacement) with probability proportional to
$\|a_{i_k}\|^2$, see \S \ref{mainresult}, {it follows that $\varepsilon_{i_k}/\|a_{i_k}\|$ is independent from 
$(b_{i_k} - \langle a_{i_k}, x_k \rangle)/\|a_{i_k}\|$.} Hence
\begin{equation} \label{expectnoise}
\zeta_k = \alpha_k^2 \sigma^2. 
\end{equation}
We remark that if
$i_k$ was chosen uniformly at random from $\{1,\ldots,m\}$ and we had
previously selected equation $i_k$, say, during iteration $j$ for $j < k$, then
the error in the $i_k$-th equation $b_{i_k} - \langle a_{i_k}, x_k \rangle$ may depend (or even be determined) by
 $\varepsilon_{i_k}$; thus the assumption that rows are drawn without replacement is necessary for this 
 term to vanish.
We use the same estimate for $\eta_k$ as in \cite{StrohmerVershynin2009}. In particular, by
\cite[eq. 7]{StrohmerVershynin2009} we have
\begin{equation} \label{etalower}
\eta_k =
\mathbb{E}_{i_0,\ldots,i_{k-1}}
\left| \left\langle \frac{x_k - x,
}{\|x - x_k\|}, \frac{a_{i_k}}{\|a_{i_k}\|} \right\rangle \right|^2
\ge \kappa(A_k)^{-2} \ge \eta,
\end{equation}
where the final inequality follows by assumption \eqref{etacond} in the statement of Theorem \ref{thm1}.

\begin{remark}[Comparison to case $\alpha_k = 1$]
Note that when $\alpha_k = 1$, the linear term {of \eqref{eq:def_zeta_k_new} satisfies}
$$
-2  \alpha_k(1 - \alpha_k)
\mathbb{E}_{i_0,\ldots,i_{k-1}} \frac{\varepsilon_{i_k}}{\|a_{i_k}\|}
\frac{b_{i_k} - \langle a_{i_k}, x_k \rangle}{\|a_{i_k}\|} = 0 
$$
 because $(1 - \alpha_k) = 0$, which geometrically is the result of an
orthogonality relation, which holds regardless of the structure of the noise.
Here we consider the case $\alpha_k \not = 1$, and the linear term vanishes due
to
the  assumption that $\varepsilon_i$ are independent mean zero random variables.
\end{remark}

\subsection{Step 4: optimal learning rate with respect to upper bound}
\label{step4}
In this section we derive
the optimal learning rate $\alpha_k$ with respect to an upper bound on the expected error.
In Remark \ref{sigma0}, we already considered the case $\sigma=0$, and found that the optimal
learning rate is to set $\alpha_k =1$ for all $k=0,1,\ldots$ regardless of the other parameters.
Thus, in this section we assume that $\sigma > 0$.
By \eqref{expectn}, \eqref{expectnoise}, and \eqref{etalower} we have
\begin{equation} \label{iterateeq}
\mathbb{E}_{i_0,\ldots,i_{k-1}}
\|x_{k+1} - x\|^2 \le \Big(1 - (2 \alpha_k - \alpha_k^2) \eta
 \Big) \|x_k - x\|^2
+ 
\alpha_k^2 \sigma^2.
\end{equation}
Iterating this estimate and taking a full expectation gives
$$
\mathbb{E}
\|x_{k+1} - x\|^2 \le g(k,\alpha) \sigma^2,
$$
for
$$
g(k,\alpha) :=
\prod_{j=0}^k
\Big(1 - (2 \alpha_j - \alpha_j^2) \eta
 \Big) \frac{\|x_0 - x\|^2}{\sigma^2}
+ 
\sum_{j=0}^k 
\alpha_j^2 
\prod_{i=j+1}^k
\Big(1 - (2 \alpha_i - \alpha_i^2) \eta
 \Big) ,
$$
where $\alpha = \{\alpha_j\}_{j=0}^{N}$ and
we use the convention that empty products are equal to $1$. {
Note that since the sums and products defining $g(k,\alpha)$ range up to $k$, it follows that $g(k,\alpha)$ does not depend on $\alpha_j$ for $j \ge k+1$.}
This upper bound $g(k,\alpha)$ satisfies the recurrence relation
$$
g(k,\alpha) = \Big(1 - (2 \alpha_k - \alpha_k^2) \eta\Big) g(k-1,\alpha) + \alpha_k^2 ,
$$
where $g(k-1,\alpha)$ does not depend on $\alpha_k$. Setting the partial derivative $\partial_{\alpha_k} g$ of $g$ with respect to $\alpha_k$ equal to zero,
and solving for $\alpha_k$ gives
\begin{equation} \label{mint}
\alpha_k = \frac{ \eta g(k-1,\alpha) }{ \eta g(k-1,\alpha) + 1}.
\end{equation}
Since $\partial^2_{\alpha_k} g(k,\alpha) =  2 \eta g(k-1,\alpha) + 2 \sigma^2 > 0$  the value of $\alpha_k$ defined by \eqref{mint} does indeed minimize $g(k,\alpha)$ with respect to $\alpha_k$.
It is straightforward to verify that this argument can be iterated to conclude that
the values of $\alpha_0, \alpha_1, \ldots $ that minimize $g(k,\alpha)$ satisfy the recurrence relation:
\begin{equation}
\begin{split}
\beta_0 &:= \|x - x_0\|^2 /\sigma^2,\\
\alpha_k &= \frac{\eta \beta_k}{\eta \beta_k + 1}
\quad \text{and} \quad
\beta_{k+1} =  \Big(1 - (2\alpha_k -
\alpha_k^2) \eta \Big) \beta_k + \alpha_k^2, \label{betaaa}
\end{split}
\end{equation}
for $k =0,1,\ldots$.
Note that we can simplify \eqref{betaaa} by observing that
\begin{equation*}
\begin{split}
    -(2\alpha_k -\alpha_k^2) \eta \beta_k + \alpha_k^2 &= \frac{-\left(2\eta \beta_k(\eta\beta_k +1) - \eta^2\beta_k^2\right)\eta\beta_k +\eta^2\beta_k^2}{(\eta \beta_k +1)^2} \\
    &= \frac{-\eta^2\beta_k^2}{\eta\beta_k + 1} = -\eta\beta_k \alpha_k.
    \end{split}
\end{equation*}
In summary, we can compute the optimal learning rate $\alpha_k$ with respect to the upper bound $g(k,\alpha)$ 
on the expected error $\mathbb{E}\|x - x_k\|^2$ as follows: if $\sigma^2 =
0$, then $\alpha_k = 1$ for all $k$. Otherwise, we define
\begin{equation}
\begin{split}
 \beta_0 &= \|x - x_0\|^2/\sigma^2,\\
\label{eq::beta}
\alpha_k &= \frac{ \eta \beta_k  }{ \eta \beta_k  + 1 }, \quad \text{and} \quad
    \beta_{k+1} = \beta_k(1-\eta\alpha_k),
    \end{split}
\end{equation}
for $k =0,1,\ldots$. In the following section we study the connection between
this recursive formula and a differential equation.

\subsection{Step 5: relation to differential equation}\label{step5}
In this section, we derive a closed form upper bound for $\beta_k$.
It follows from \eqref{eq::beta} that
$$
\frac{\eta \beta_{k+1} - \eta \beta_k}{\eta} =  - \eta \beta_k  \frac{\eta \beta_k}{\eta \beta_k + 1}.
$$
Making the substitution $u_k := \eta \beta_k$ gives the finite difference equation 
$$
\frac{u_{k+1} - u_k}{\eta} = -\frac{u_k^2}{u_k+1},
$$
which can be interpreted as one step of the forward Euler method (with step size
$\eta$) for the ordinary differential equation
\begin{equation}\label{eq::def_x}
\dot{u} = - \frac{u^2}{u+1},
\end{equation}
where $u = u(t)$ and $\dot{u} = du/dt$. 
%%We next write down an explicit expression for an approximation of the $\beta_k$ defined through the recursion relation in~\eqref{eq::beta}. This is done through the solution of a closely related differential equation for the auxiliary quantity $x_k = \eta \beta_k$.
It is straightforward to verify that the solution of this differential equation is 
\begin{equation}\label{eq::sol_x}
u(t) = \frac{1}{W\left(e^{t+c} \right)},
\end{equation}
where $W$ is the Lambert-$W$ function (the inverse of the function $x \mapsto
xe^x$) and $c$ is determined as the initial condition; in particular, if $u(0) =
u_0$, then
\begin{equation} \label{init}
c = \frac{1}{u_0} - \ln(u_0).
\end{equation}
We claim that $u$ is a convex function when $u(0) \ge 0$. It suffices to check that 
$\ddot{u} \ge 0$. Direct calculation gives
\begin{equation} \label{ypp}
\ddot{u} = \frac{u^4 + 2u^3}{(u+1)^3}.
\end{equation}
Observe that
$u$ cannot change sign because $\dot{u} = 0$ when $u = 0$. Thus, \eqref{ypp} is always nonnegative when $u(0) \ge 0$ as
was to be shown. Since the forward Euler method
is a lower bound for convex functions, it follows from \eqref{eq::sol_x} and \eqref{init} that
$$
\beta_k \le \frac{1}{\eta W(e^{\eta k+c})},
\quad \text{for} \quad 
 c := \frac{1}{\eta \beta_0} - \ln( \eta \beta_0).
$$
Thus if we set
\begin{equation}\label{eq::fdefinproof}
f(k) :=
 \frac{\sigma^2}{\eta
%%% W(e^{\eta k+c})},
W(c e^{\eta k})},
\quad \text{for} \quad
%%% c := \frac{\sigma^2 }{\eta \|x - x_0\|^2} - \ln( \eta \|x - x_0\|^2/\sigma^2),
c := \frac{\sigma^2 }{\eta \|x - x_0\|^2} \exp \left( \frac{\sigma^2 }{\eta \|x - x_0\|^2} \right),
\end{equation}
it follows that 
\begin{equation} \label{lasteq}
\mathbb{E} \|x - x_k\| \le \sigma^2 \beta_k \le \frac{\sigma^2}{\eta} u(\eta k) = f(k).
\end{equation}
This completes the proof of Theorem \ref{thm1}.

\section{Proof of Corollaries} \label{proofcor}

\subsection{Proof of Corollary \ref{cor1}} \label{proofcor1}

By the definition of $f(k)$ and $c$, see \eqref{eq::fdefinproof}, we have
$$
f(k) = \frac{\sigma^2}{\eta W\left(
\frac{\sigma^2}{\eta \|x - x_0\|^2} e^{ \eta k}  e^{\sigma^2/(\eta \|x - x_0\|^2)} \right)}.
$$
The Lambert-$W$ function has the Taylor series
$$
W(x) = \sum_{n=1}^\infty (-1)^{n-1} \frac{n^{n-2}}{(n-1)!} x^n,
\quad \text{for} \quad |x|<\frac{1}{e},
$$
see for example \cite[eq. 4.13.5]{NIST:DLMF}; in particular, $W(x) = x +
\mathcal{O}(x^2)$ as $x \rightarrow 0$. Thus,
$$
f(k) = \frac{\sigma^2}{\eta} \frac{1}{
\frac{\sigma^2}{\eta \|x - x_0\|^2} e^{ \eta k}  e^{\sigma^2/(\eta \|x - x_0\|^2)} + \mathcal{O}\left( \frac{\sigma^4}{\|x - x_0\|^4} \right)}.
$$
Canceling terms and using the fact that $e^{\sigma^2/(\eta \|x - x_0\|^2)} = 1 + \mathcal{O}(\sigma^2/\|x-x_0\|^2)$ gives
$$
f(k) = e^{-\eta k} \|x - x_0\|^2 \left( 1 + \mathcal{O}\left(\frac{\sigma^2}{\|x-x_0\|^2} \right) \right), \quad \text{as} \quad \sigma \rightarrow 0,
$$
where the constant in the big-$\mathcal{O}$ notation depends on $\eta$ and $k$, as was to be shown.

\subsection{Proof of Corollary \ref{cor2}} \label{proofcor2}
By the definition of $f(k)$ and $c$, see \eqref{eq::fdefinproof}, we have
$$
f(k) = \frac{\sigma^2}{\eta} \frac{1}{W\left(
\frac{\sigma^2}{\eta \|x - x_0\|^2} e^{ \eta k}  e^{\sigma^2/(\eta \|x - x_0\|^2)} \right)}.
$$
The Lambert-$W$ function has asymptotic expansion 
$$
W(e^\xi) = \xi - \ln(\xi) + \mathcal{O} \left( \frac{\ln \xi}{\xi} \right), \quad \text{as} \quad \xi \rightarrow +\infty,
$$
see \cite[eq. 4.13.10]{NIST:DLMF}. It follows that
$$
f(k) = \frac{\sigma^2}{\eta^2 k} \left( 1 + \mathcal{O}\left( \frac{\ln k}{k} \right) \right),
\quad \text{as} \quad k \rightarrow \infty,
$$
where the constant in the big-$\mathcal{O}$ notation depends on $\eta$, $\|x-x_0\|^2$ and $\sigma^2$, as was to be shown.

\subsection{Proof of Corollary \ref{cor3}}\label{proofcor3}
The proof of Corollary \ref{cor3} is immediate from the proof of Theorem \ref{thm1}; in particular, see \eqref{etalower}.

\subsection{Proof of Corollary \ref{cor4}} \label{proofcor4}
First we argue why \eqref{cond} holds with equality when it holds with $\eta = 1/n$.
Let $z_1$ be an arbitrary unit vector and complete it to an orthonormal basis $\{z_1, z_2, \ldots , z_n\}$. By assumption, the expected squared magnitudes of the coefficients of $a/\|a\|$ in this basis satisfy
\begin{equation}\label{eq::etacondproof4}
\mathbb{E}
\left| \left\langle z_i, \frac{a}{\|a\|} \right\rangle \right|^2 \geq \eta = \frac{1}{n}.
\end{equation}
The sum of the squares of the coefficients of $a/\|a\|$ in any orthonormal basis is
equal to $1$. It follows that \eqref{eq::etacondproof4} holds with equality for
each $i=1, \ldots , n$, and in particular for $z_1$. Since $z_1$ was arbitrary,
\eqref{eq::etacondproof4} holds with equality for arbitrary unit vectors.
If \eqref{etalower}
in the proof of Theorem \ref{thm1} holds with equality, then it is
straightforward to verify that the remainder of the proof of Theorem \ref{thm1}
also carries through with equality, so the bound in
Theorem \ref{thm1} is sharp in this case, which concludes the proof.

\section{Discussion} \label{discuss}
In this paper, we have presented a randomized Kaczmarz algorithm
with a scheduled learning rate for
solving $A x \approx b + \varepsilon$,
where $A x =b$ is a consistent linear system and $\varepsilon$ has independent mean zero random entries. 
When we start with $x_0 = \vec{0}$, the scheduled learning rate $\alpha_k$ defined by \eqref{optlearn} depends on two parameters:
\begin{itemize}
\item the signal-to-noise ratio $\|x\|^2/\sigma^2$, and
\item the condition number parameter $\eta$.
\end{itemize}
This learning rate optimizes the error bound of Theorem \ref{thm1} which is sharp in certain cases, see Corollary \ref{cor4}.
There are many extensions of the randomized Kaczmarz algorithm of \cite{StrohmerVershynin2009} such as
\cite{
eldar2011acceleration,
Liu2015,
Ma2015,
Necoara2019,
Needell2015b,
needell2014paved,
Petra2016,
Tan2018,
Zouzias2013}, which could be considered in the context of our model and analysis.
In particular, it would be interesting to consider the  block methods of \cite{
Necoara2019,
Needell2015b,
needell2014paved,
slagel2019row,
slagel2019sampled}. In the context of machine learning, blocks correspond to batches which are critical to the performance of stochastic gradient descent in applications. 
In the same direction, connections to adaptive learning rates such as 
Adadelta \cite{Zeiler2012} and ADAM \cite{Kingma2014} would also be interesting to consider.

In practice, optimization algorithms are run with epochs. In the context of our method, after looping over the data once, we can set $x_0$ using our final iterate and loop over the data again. Formally, the statement of Theorem \ref{thm1} no longer holds, but practically, the iteration error may continue to decrease in some cases. In particular, if the iteration error has not reached a ball of radius $\sigma/ \sqrt{\eta}$ around the solution (see \eqref{needell}), then practically speaking, \eqref{expectnoise} might still approximately hold, and the result of the theorem might still carry through. This could potentially be studied with a more detailed analysis.

\subsection*{Acknowledgements}
We are grateful to Marc Gilles for many helpful comments. 
{
We also thank the anonymous reviewers for their insightful comments, which greatly improved the exposition of the results.}

\end{document}